# Dimensional analysis in forest mensuration
by
Tim Davis[1] & Huw Williams[2]

## Introduction

The "`trees`" dataset, which is described in Appendix 1, is a well-known regression teaching example that (to the best of our knowledge) was originally introduced to the statistical community in the first edition of the Minitab Student Handbook in 1976[3]. The `trees` dataset (which relates to 31 Black Cherry trees from the Allegheny National Forest in Pennsylvania) has been posted in many locations, the R dataset package perhaps being the most well-known[4]. As well as being used in teaching, the dataset has been the subject of published research[5]. We believe that this dataset can teach us a lot about how to apply statistics well in solving problems by making the most of the scientific and engineering context of data, and this approach is the focus of this article. In this spirit, we managed to obtain a copy of the 1976 first edition of the Minitab handbook, which in turn led us to the original source of the data[6] : a 1953 textbook by the forester AH Meyer[7]). Henceforth we refer to this data as Meyer's cherry tree data, rather than the Minitab cherry tree data.

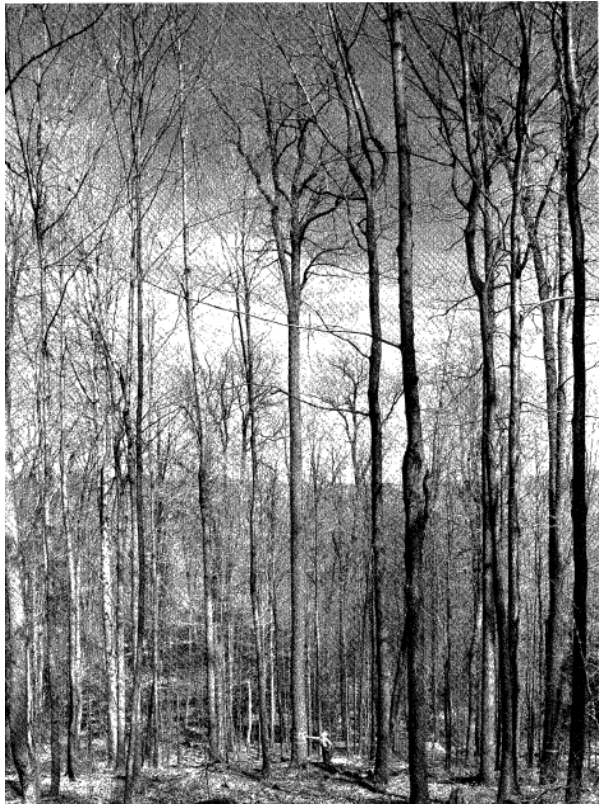

<u>**Figure 1:**</u> A forest of Black Cherry Trees[8] The tallest cherry tree ever recorded in the Allegheny National Forest (as of 1960) was 129 feet and a tree takes about 60 years to grow to 80 feet.

---

[1] Corresponding author timdavis consulting ltd & Dept of Statistics, University of Warwick, tim@timdavis.co.uk
[2] Huw Williams consulting & Dept of Engineering, University of Birmingham
[3] *MINITAB Student handbook*, by TA Ryan, B. Joiner & BF Ryan 1976. Duxberry Press (first edition)
[4] Given that this data set is provided in R suggests that this historical data is still used for teaching and research, and so ongoing investigation of approaches to the analysis of this data is still worthwhile.
[5] see Atkinson, A.C., 1994. "Transforming both sides of a tree". *The American Statistician*, 48(4), pp.307-313.
[6] The reference to the cherry tree data has been omitted from subsequent editions.
[7] Meyer, H.A., 1953. *Forest mensuration.* Penns Valley Publishers Inc.
[8] Taken from "Black Cherry Provenances for Planting in Northwestern Pennsylvania" by Russell L Walters (The Northeastern Experiment Station part of the Forest Service in the US department of agriculture (1985))

None of the statistics literature says anything about the trees themselves, which are of the cherry family. The Black Cherry is the largest of the native cherry trees of the United States and is cultivated for its timber (lumber) rather than for its wood pulp or its fruit. Meyer's text confirms that it is to this type of cherry tree that the data refers. These trees are tall with a long straight trunk and a crown of leafy branches that forms part of the forest canopy as shown in Figure 1.

Determining trunk volumes is known in forestry by the term "forest mensuration"[9]. There is a vast amount of literature on forest mensuration, and in addition to Meyer's text the interested reader can find an approachable summary here[10]. The first thing to note is that tree volume is the volume of the trunk because neither the crown of the tree, nor the stump is used[11]. The trunk from the stump to the crown is referred to as the stem of the tree. (there is also a marketable wood volume that considers the sawmill waste when the tree is milled into boards): we can't find anything in the literature that confirms this, but the forestry literature indicates that it is the case. The procedure for measuring tree volume after the tree is felled is to cut the trunk (stem) into 16ft logs. The volume of each log is then estimated by multiplying the cross-sectional area at the mid-point (estimated by either measuring the diameter at the mid-point or using the mean of the diameters of the log's ends) by its length. The resulting volumes are then summed to give the overall volume (see Appendix 3).

For the 31 trees under discussion the diameter of the trunk is measured at breast height (4.5ft) either with a girthing tape or a calliper, and then the height is estimated via basic trigonometry with a clinometer[12] (see Figure 2), the tree is then felled, and the volume determined using the procedure already described. The purpose of our analysis (and others) is to estimate this volume via a statistical model based on the measurements of diameter and height, so that the model can be used to predict the volume of a tree stem *before* it is felled (for obvious reasons), as part of ongoing forestry management, including the sale of timber to produce furniture and other goods.

---

[9] Forest Mensuration is that branch of forestry that deals with the determination of measurements (e.g., diameter, height, volume etc.), form, age, and increment of single trees, stands or whole woods, either standing or after felling.
[10] Garcia, O., 1995. Notes on forest mensuration I. Statics. Institute of Forest Management Austral University of Chile, pp.7-27.
[11] Predicting the biomass of the tree would need to include the crown and the stump, but that is not of concern here.
[12] The cherry tree data were obtained in the early 1950's, before the advent of laser devices.

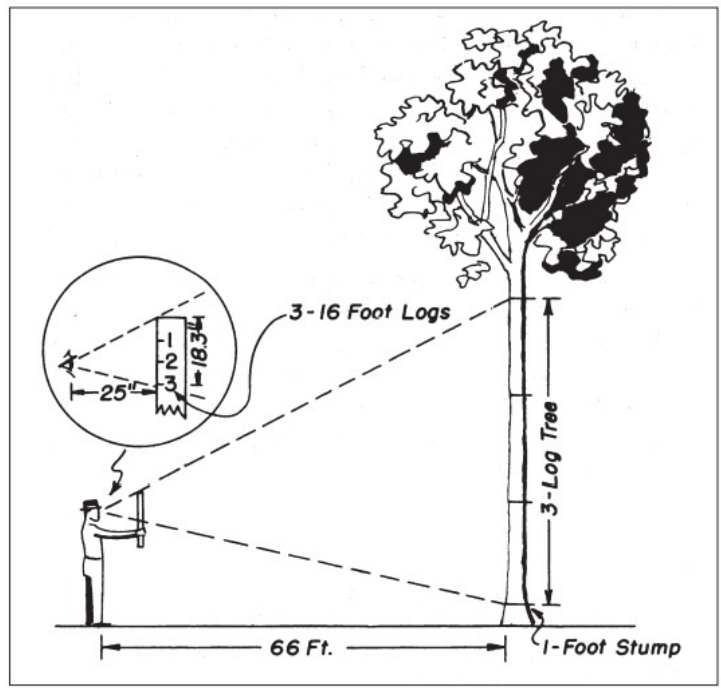

Figure 2: a schematic for the use of a clinometer to measure tree height. We assume that this is how the heights were measured for this data given that it was originally published in 1953. These days laser-based instrumentation and smart phone apps are used, but the principle is the same. Note that this diagram focuses on graduating the height in multiples of 16-foot log lengths, which is the standard log length in forestry.

The data was originally presented as a study in statistical modelling with the goal of predicting the volume of standing cherry trees before they are felled. There are several discussions in the statistical literature of fitting these statistical models which we will refer to as we discuss our analysis.

## Modelling the cherry tree data

Meyer's Cherry tree data (see Appendix 1) is illustrated in Figure 3.

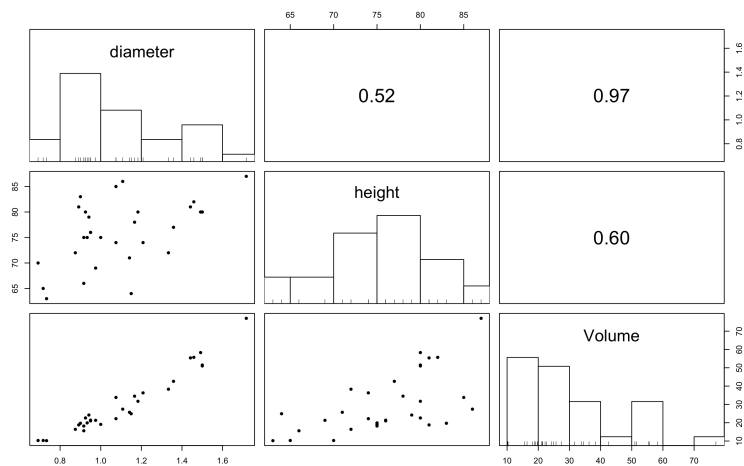

Figure 3: a plot of Meyer's cherry tree data. In this plot we have expressed Diameter in feet for consistency of units with height and volume. The numbers in the upper triangle panels are the correlation coefficients.

If we conjecture (based on previous studies of this data[13]) that the trunk of a cherry tree is either cylinder or cone shaped (or something in between), we might, because of considerations of dimensions, expect to see a relationship of the form $V \propto hd^2$, where the constant of proportionality would be $\frac{\pi}{12}$ for a cone and $\frac{\pi}{4}$ for a cylinder[14]. Prior to fitting a model and given we are dealing with trees we might be tempted to take logs which leads to a linear model[15] of the form

$$\log(V) = \beta_0 + \beta_1 \log(d) + \beta_2 \log(h) + \epsilon,$$

where $\epsilon$ captures the errors (departures from reality) for this model[16], if our conjectured models are close to being correct, we might expect that the coefficients $\beta_1 \approx 2$, and $\beta_2 \approx 1$. We will refer to this model as the log-regression model.

Fitting this model by least squares yields the following estimates denoted by $\hat{\beta}_0$ for $\beta_0$ etc. $\hat{\beta}_0 = -1.705_{(0.8819)}; \hat{\beta}_1 = 1.98_{(0.0750)}; \hat{\beta}_2 = 1.117_{(0.2044)}$. The subscripted numbers in parentheses are the standard errors of the estimates. Note that both $\hat{\beta}_1$ and $\hat{\beta}_2$ are within a single standard error of 2 and 1 respectively, which at first sight might lend weight to the cylinder/cone type model. Now $\log(\frac{\pi}{12}) = -1.3408$, while $\log(\frac{\pi}{4}) = -0.2416$, and since $\hat{\beta}_0$ is within a single standard error of $\log(\frac{\pi}{12})$, it might be tempting to suggest that the log-regression model is suggesting that a cone is the appropriate model. But things are not that straightforward. In a previous treatment Fairley[17] has pointed out that $(\hat{\beta}_0, \hat{\beta}_1, \hat{\beta}_2)$ are not independent[18] and so individual statistical assessments of each $\hat{\beta}$ as if they were, is likely to be misleading. In fact, the correlation between $\hat{\beta}_0$ and $\hat{\beta}_2$ is -0.9998 and around +0.5 for the other two pairings.

These correlations need to be considered when constructing uncertainty intervals for the model coefficients, and this requires the construction of a confidence ellipsoid (rather than just a single interval) which considers the correlation between $\hat{\beta}_0$ and $\hat{\beta}_2$. For ease of exposition, we have reproduced in Figure 2 the confidence ellipsoid given by Figure 1 of Fairley's paper.

---

[13] e.g. see AC Atkinson (1985) *Plots Transformation and regression* Oxford Science publications, and Cook & Weisberg (1982) *Residuals and influence in regression* Chapman & Hall

[14] Given the strong correlation between $d$ and $V$ ($\rho_{dV} = 0.97$, see Figure 1), it might be tempting to consider a model for $V$ just in terms of $d$ since $d$ is very easy to measure and would be simple for the forester to use. Dimensional considerations would suggest that we could try $V = f(d^3)$ for some $f(.)$ to be determined, but it turns out that this model has high prediction error and a plot of $V$ versus $d^3$ does not pass through the origin as one might expect, so we do not consider it further., and concentrate henceforth on models that depend on $h$ as well as $d$.

[15] Note linear here means linear with respect to the $\beta's$

[16] $\epsilon$ is specified to have zero mean and variance $\sigma_\epsilon^2$, and in the context of least square regression is taken as Gaussian; $\sigma_\epsilon^2$ can be estimated once a model has been fitted.

[17] David Fairley – "Cherry trees with cones?", *The American Statistician*, May 1986, Vol 40. #2

[18] The $\hat{\beta}'s$ are functions of the data so they inherit a variance structure transmitted from the data, which will be available from statistical software.

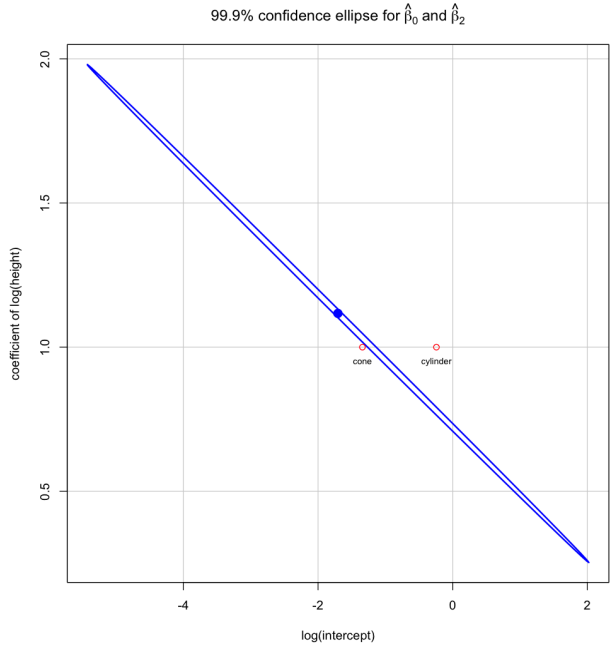

**Figure 2:** A reproduction of the 99.9% confidence ellipsoid for the regression model, as given by Fairley. Note that both the cylinder and cone model fall outside the ellipsoid. The log-regression model is of course at the centre of the ellipsoid, illustrated by the blue dot.

Figure 2 shows that, even though both models are dimensionally consistent, we can reject the cylinder and cone models as not being consistent with (in the sense of implied by) the log-regression model.

We can also reject both the cylinder and cone models based on considerations of the physics and geometry involved. Firstly, we reject the cone model because we can imagine no tree whose trunk diameter approaches zero directly below the crown of the tree, and we reject the cylinder model because trunk diameters taper at higher points on the trunk, a consequence of natural selection requiring the evolution of stable structures, to deal with the effects of gravity.

### Dimensional Analysis

So how might we proceed to develop a dimensionally consistent model that will be useful in predicting $V$ as a function of $h$ and $d$? Since we have three variables ($V, d, \& h$) with just one unit (feet) Buckingham's Pi theorem (see Appendix 2) states that we can reduce the problem to $3 - 1 = 2$ dimensionless variables. The choice of these dimensionless variables is not unique (as with the introductory toilet roll example[19] discussed in Appendix 2). One choice is to pick $\pi_0 = \frac{V}{h^3}$, and $\pi_1 = \frac{d^2}{h^2}$, and another is $\pi_0 = \frac{V}{d^3}$, and $\pi_1 = \frac{h}{d}$. Figure 4 illustrates a total of four obvious possibilities[20]

---

[19] It is worth pointing out that black cherry trees are not used in the production of toilet paper!
[20] There are more possibilities but these four are the ones that seem at first sight to have the most utility.

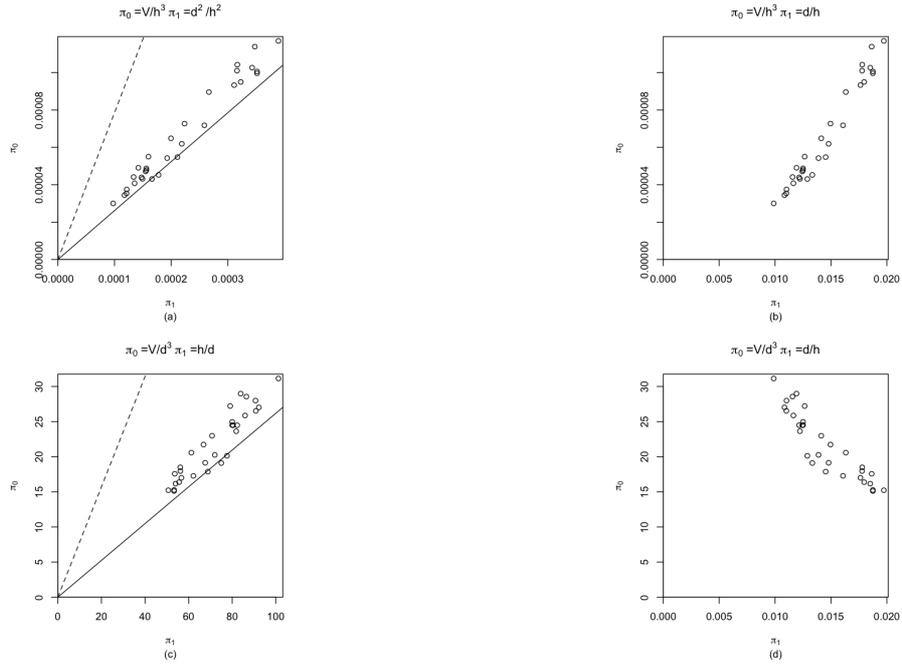

**Figure 4:** Plots of the cherry tree data with various choices for the two dimensionless variables. For formulations (a) and (b) the dotted line represents a cylinder and the solid line a cone. This way of representing the data clearly demonstrates that neither a cylinder nor a cone is appropriate for this data. For formulation (b) the data do not point to the origin, and neither a cylinder nor a cone can be represented by formulation (d).

The task now is to fit a linear model $\pi_0 = f(\pi_1)$, which we will call the DA model. Monomial models of the type $\pi_0 = \gamma_0 \pi_1^{\gamma_1}$ are preferred if only because *all* derived units take this form with integer powers[21]. For monomial models fitted to dimensionless variables, the $\gamma$ coefficients representing exponents will consequently be contained within the set of rational numbers expressed as simple fractions.

Inspection of Figure 4 would suggest that for the two representations in (a) and (c) are of the monomial form $\pi_0 = \gamma_0 \pi_1$ (i.e., $\gamma_1 = 1$)[22].

So, we fit the model $\pi_0 = \gamma_0 \pi_1$ directly (i.e., a straight line through the origin with $\pi_1$ on the x-axis and $\pi_0$ on the y-axis).

If we do this using straightforward least squares regression, we find that for (a) $\hat{\gamma}_0 = 0.302355_{(0.003893)}$, and for (c) $\hat{\gamma}_0 = 0.30270_{(0.00423)}$.

We slightly favour the model given by the DA framework (a) if only because the standard error (in subscripted parentheses) of the $\gamma_0$ coefficient is smaller.

So, our fitted DA model is $\pi_0 = 0.302 \pi_1$, or $\frac{V}{h^3} = 0.302 \frac{d^2}{h^2}$, hence $V = 0.302 h d^2$.

Note that unlike the log-regression model, the dimensional analysis approach *forces* the exponents of $h$ and $d$ to be 1 and 2 respectively.

---

[21] See Szirtes T (2007) *Applied Dimensional Analysis and modelling* Elsevier. For example, the derived unit Volt has units kgm²s⁻³A⁻¹.
[22] This can be confirmed by fitting the linear model $\log(\pi_0) = \gamma_0^* + \gamma_1 \log(\pi_1)$, whereupon $\gamma_1$ is estimated as 1 to two significant figures for both representations Figure 3(a) and (c).

If we add this model to the confidence ellipsoid previously illustrated in Figure 1, we find that this model fits inside the confidence region (see Figure 5).

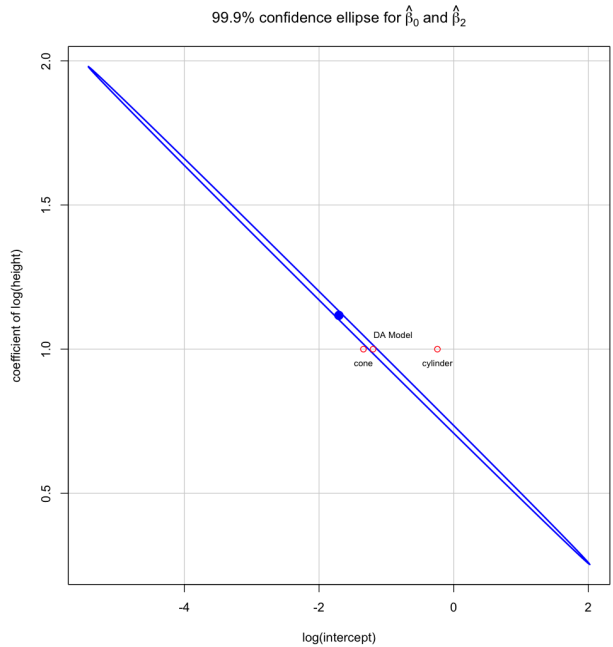

**Figure 5:** The DA model added to the confidence ellipsoid, which shows that unlike the cone and cylinder the DA model is consistent with the log-regression model.

Interestingly in an earlier DA treatment for the cherry tree data, Vignaux & Scott[23] choose the representations (b) and (d) for their initial DA models, but these result in preliminary models $V = \gamma_{0b}h^3 + \gamma_{1b}dh^2$ for (b) and $V = d^3\left(\gamma_{0d} + \gamma_{1d}\frac{d}{h}\right)$ for (d); and although these models are dimensionally consistent, as they must be, they are not monomial, so we reject these choices on this basis without further investigation.

Because of the agreement of the DA model with the log-regression model illustrated in Figure 5, regression diagnostics[24] will likely be similar for our DA model and for the log-regression model, so we won't overload our analysis with discussing model diagnostics here, but we invite the interested reader to check this for themselves (we ourselves found nothing untoward).

---

[23] GA Vignaux & JL Scott "Simplifying regression models using dimensional analysis", *Austral. & New Zealand J. Statist.* 41(1), 1999, 31–41
[24] i.e., plots of residuals (the differences from the actual values of $V$ from the predicted values determined from the model) in various forms.

## Model interpretation

The question now arises as to how we should interpret our model? It is of course dimensionally consistent, but as shown in Figure 4 it is neither a cone nor a cylinder. Tree trunks though, are tapered. Suppose the diameter of the trunk just under the crown is denoted by $\lambda d$, where $\lambda \in (0,1)$. A tapered cylinder is called a frustum of a cone. The volume of a frustum is given by $V = \frac{\pi}{12}(1+\lambda+\lambda^2)d^2h$. The volume of the cylinder and cone can be recovered by setting $\lambda = 1$, and $\lambda = 0$ respectively. For our model we can equate $\frac{\pi}{12}(1 + \lambda + \lambda^2)$ with our coefficient 0.302, re-arranging to solve the quadratic gives $\hat{\lambda} = 0.13526$, where, as usual the 'hat' denotes an estimate, and which can be interpreted as the average $\lambda$ across the 31 trees. Foresters define taper as $t = -\left(\frac{d-\lambda d}{h}\right) = (\lambda - 1)\frac{d}{h}$; if we calculate this taper value for the "average tree" (a tree with a diameter of 1.1 feet, and a height of 76 feet) we find that $t = -0.0126$, which is not too far from the figure of -0.0098, given in a key forestry reference[25] where this figure is derived by measuring taper directly on a sample of 79 trees. Ideally, we would like to calculate taper directly for the cherry trees, but this would require a measurement of the stem diameter just below the crown for each tree, and this data is not available[26].

## Transmitted variation

Finally, we recognise that given the way the trees are measured prior to felling it is of interest to study how errors in measuring $d$ and $h$ transmit to estimating $V$[27] via the well-known variance transmission formula:

$$\sigma_V^2 = \left(\frac{\partial V}{\partial d}\right)^2 \sigma_d^2 + \left(\frac{\partial V}{\partial h}\right)^2 \sigma_h^2 + \left|\frac{\partial V}{\partial d}\right| \cdot \left|\frac{\partial V}{\partial h}\right| \rho_{dh} \sigma_d \sigma_h$$

The last term is required because $d$ and $h$ are correlated ($\rho_{dh} = 0.52$). Note that the units of $\sigma_V^2$ are ft$^6$ and the units of $\sigma_V$ ft$^3$. This equation is the cosine rule for variances. However, to investigate $\sigma_V^2$ given by the above equation we need an estimate of both $\sigma_d$ and $\sigma_h$., which would require several measurements of $d$ an $h$ on the *same* tree by *different* foresters, conducted via a gauge R&R (repeatability and reproducibility) study[28] which we don't have. Meyer discusses the accuracy of measuring $d$ and $h$ and although he doesn't provide values for $\sigma_d$ or $\sigma_h$, he does imply that measuring $h$ is inaccurate due to trees not always growing upright, not always being clearly visible in a forest, not always growing on level ground, etc., while in contrast measuring $d$ is relatively accurate using a girthing tape or a trunk calliper, whereby inaccuracies creep in due to trunks not being perfectly circular, bark not being smooth[29] etc. For us to illustrate the variance transmission formula, we will assume that we can measure $h$ within $\pm$5%, and $d$ within $\pm$1%, giving implied coefficients of variation[30] for $h$ and $d$ as 0.0408 and 0.0082 respectively.

---

[25] See Table 1 in "Simple Taper: taper equations for the field forester" by David R. Larsen, in *Proceedings of the 20th Central Hardwood Forest Conference.*

[26] Although we speculate that it was measured during sectioning the stem into logs, so it is regrettable that this aspect of the data is not available.

[27] This would be important for the forester to know for minimizing risk in a contractual sale.

[28] For example, see https://en.wikipedia.org/wiki/ANOVA_gauge_R%26R.

[29] Meyer indicates (pg.54) that callipers tend to underestimate the diameter if the arms of the device are not exactly parallel, while a girthing tape tends to overestimate the diameter. Since it is not clear which method was used for the cherry tree data, we adopt a $\pm$% treatment for the error. He discusses errors in measuring height on page 101 of his text.

[30] The coefficient of variation is $c_v = \frac{\sigma}{\mu}$, which if constant as assumed here implies that $\sigma$ increases in proportion to the mean, which is a consequence of using $\pm$% error formulation. So $\sigma = c_v \mu$, which can be substituted into the variance transmission formula.

Figure 6 shows the contours of the transmitted variation from $d$ and $h$ to $V$ using these guesses for measurement accuracy. Note that the transmitted variance is far more sensitive to errors in $d$ than $h$ *in spite of the fact that $d$ is measured with more accuracy.*

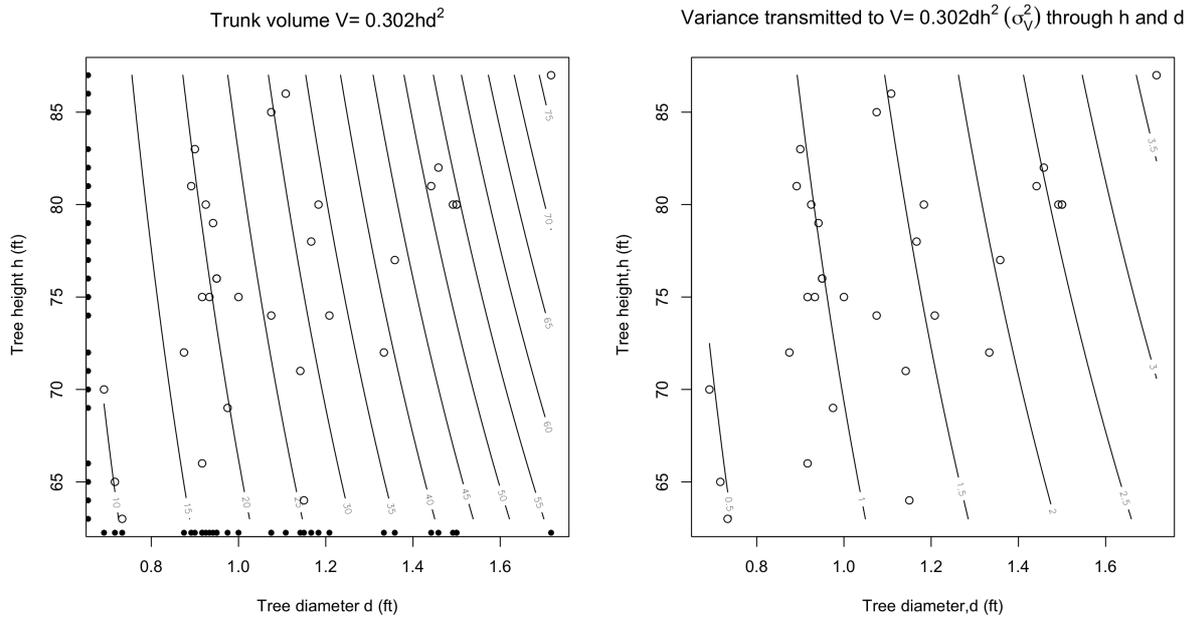

**Figure 6:** contours of $V$ (left) and $\sigma_V^2$ (right) as a function of $d$ and $h$. Assuming $h$ can be measured to an accuracy of $\pm$5% and $d$ $\pm$1%. The dimensions of the 31 trees have been superimposed, onto the contours. Note that from the right-hand plot the variance in $V$ for the largest tree due to these assumed errors in measuring $d$ and $h$ is about 3.5ft$^6$

## Concluding remarks

We have emphasised the importance of a thorough understanding of how data is collected before embarking on a statistical analysis. There is no substitute from being there while the data is collected so the data collection can be observed first hand. Of course, this is not possible for this data which is why we embarked on a thorough review of the forestry literature so that we could understand the engineering and scientific principles involved in felling trees, in that way we could imagine we were there before we attempted an analysis. We hope that our article encourages educators to take sufficient care with teaching examples so that students appreciate the importance of the context of how their subject data is collected. We call our approach *statistical engineering*.

The use of Buckingham's "Pi" theorem reduces the number of variables in the modelling phase, and this encourages model parsimony, leading to more realistic models. We are not aware of previous studies that have entertained a frustum of a cone for this data. We claim that this model is realistic because it is in close agreement with already published data on cherry tree taper published by the US forestry service.

There is a saying attributed to George Box "All models are wrong, but some are useful". Although our model predicts the recorded volumes of the cherry tree data well[31] the usefulness of our model can only really be validated by collecting a new set of data on cherry trees in the Allegheny National Forest at a time relatively close to when the data was originally obtained, which of course we are not able to do. What we can do is to try our model on another tree data set. Fortunately, a data set on 70 short leaf pine trees exists in the reference noted in footnote 3. If we fit our DA model to the pine tree data, we find that $\hat{\gamma}_0 = 0.43629 \ _{(0.00363)}$. Although this is a different species of tree altogether, the scalability of the DA model is apparent. See figure 7 for a comparison of the cherry and pine tree data.

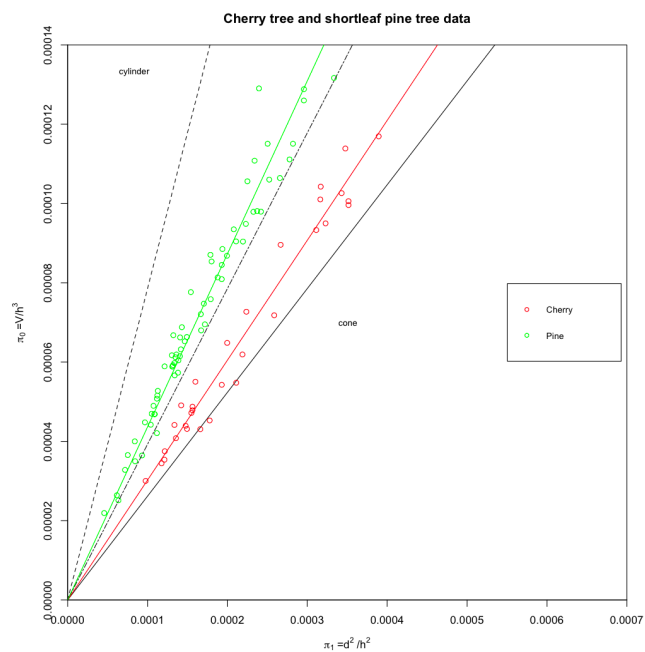

**Figure 7:** a comparison of the cherry tree data with the shortleaf pine data in terms of the dimensionless variables, with DA models superimposed over the data. This plot illustrates the scalability of dimensional analysis, in the sense that the same model (with adjusted slope) fits both species of tree.

Further investigation shows that our model reproduces published tables for predicted volumes very well with minimal error. For example, Table 6.7 in the text *Forest Mensuration*[32] has volume tables (in m$^3$) based on $d$ (graduated in 2cm increments) and $h$ (in m) on Spruce, Balsam Fir, White Pine, Maple, White Birch, and Yellow Birch trees from the Noonan Forest in Central New Brunswick, Canada, where the tabulated stem volumes have been calculated using Honer *et al's* formula[33] Note that with our DA

---

[31] The prediction error (residual) sum-of-squares is $\sum_1^{31}(V_i - 0.302 h_i d_i^2)^2$ for our DA model is 181.4 cubic feet.

[32] John A. Kershaw, Jr, Mark J. Ducey, Thomas W. Beers, Bertram Husch *Forest Mensuration* 5th edition (2017); Wiley: Blackwell see Table 6.7 page 163.

[33] See Honer, T. G., M. F. Ker, and I. S. Alemdag. 1983. *Metric timber tables for the commercial tree species of central and eastern Canada. Information Report M-X-140, Canadian Forestry Service, Maritimes Forest Research Centre.* (Available for download from the internet); These authors' volume equation is $V = \left(\frac{d^2}{c_1 + c_2/h}\right)$ (their eq. 13) where $c_1$ and $c_2$ are estimated parameters from previous mensuration studies. For Black Cherry trees these parameters are given in their

approach we do not need to convert these measurements to feet for comparison with what we have already done (we do however, express $d$ in m for consistency of units). The DA models for all these trees fall between the Pine and Cherry tree models illustrated in Figure 7, illustrating the wide applicability of our DA model. We believe modelling different tree species in this way, could lead to some interesting scientific insights, connecting our DA model parameters with geo-biological measures such as species, latitude & longitude, and growth rate. We would be delighted to correspond with current-day foresters to investigate applications and further extensions to the DA model.

There is a further data set in Meyer's text concerning 112 Hartweg Pines (Table 36, pg. 128). If readers wish to try our method on another tree dataset, that hasn't received any prior attention it is available from the corresponding author.

---

Table 3 (pg. 17) as 0.033 and 393.336 respectively, based on a study of 21 trees. Note that in this equation $d$ is expressed in inches (as measured), $h$ is in feet and $V$ is in ft³; it is not monomial and is non-linear in the parameters. It fits Meyer's cherry data remarkably well, and by implication, we suppose that our DA model would predict $V$ for the 21 cherry trees cited to a similar degree.

# Appendix 1: The cherry tree data

Diameter, height and volume of 31 felled trees of Black Cherry, Allegheny national forest, Pennsylvania

| Diameter breast height | Total Height H | Actual Volume V |
|---|---|---|
| inches | feet | cubic feet |
| 8.3 | 70.0 | 10.3 |
| 8.6 | 65.0 | 10.3 |
| 8.8 | 63.0 | 10.2 |
| 10.5 | 72.0 | 16.4 |
| 10.7 | 81.0 | 18.8 |
| 10.8 | 83.0 | 19.7 |
| 11.0 | 66.0 | 15.6 |
| 11.0 | 75.0 | 18.2 |
| 11.1 | 80.0 | 22.6 |
| 11.2 | 75.0 | 19.9 |
| 11.3 | 79.0 | 24.2 |
| 11.4 | 76.0 | 21.0 |
| 11.4 | 76.0 | 21.4 |
| 11.7 | 69.0 | 21.3 |
| 12.0 | 75.0 | 19.1 |
| 12.9 | 74.0 | 22.2 |
| 12.9 | 85.0 | 33.8 |
| 13.3 | 86.0 | 27.4 |
| 13.7 | 71.0 | 25.7 |
| 13.8 | 64.0 | 24.9 |
| 14.0 | 78.0 | 34.5 |
| 14.2 | 80.0 | 31.7 |
| 14.5 | 74.0 | 36.3 |
| 16.0 | 72.0 | 38.3 |
| 16.3 | 77.0 | 42.6 |
| 17.3 | 81.0 | 55.4 |
| 17.5 | 82.0 | 55.7 |
| 17.9 | 80.0 | 58.3 |
| 18.0 | 80.0 | 51.5 |
| 18.0 | 80.0 | 51.0 |
| 20.6 | 87.0 | 77.0 |

Data on 31 cherry trees collected in the Allegheny National Forest Pennsylvania, as the data appears in Table 32 in the textbook *Forest Mensuration* by AH Meyer (pg. 118). Not shown here are two extra variables tabulated by

Appendix 1 (cont:)

Meyer, *Balanced height* (which is a prediction of tree height based on diameter), and *Adjusted Volume* (which uses balanced height in lieu of measured height to estimate volume) These are derived quantities, and not of interest here.

The preface to the book reads" The development of new techniques in forest mensuration through the application of statistical methods…has greatly increased the scope of this basic field of forestry". While discussing the cherry data Meyer initially suggested a model only involving the diameter at breast height (which foresters typically refer to as DBH, but which we will refer to as simply $d$, with Volume as $V$ and height as $h$) of the form $V = kd^b$, which surprisingly given his comments in the preface he took no further. Dimensional consistency would require $b = 3$.

Conventionally breast height is taken at 4.5ft above the ground, thus avoiding the tree stump. Foresters usually refer to the trunk above the stump as the *stem* of the tree. All our analysis converts diameter to feet for consistency in units since height is recorded in feet and volume in cubic feet. We have resisted the temptation to convert these data to SI units (i.e., metres and cubic metres) to facilitate comparisons of our treatment of this data to what has been published previously.

# Appendix 2: Buckingham's Pi Theorem

Maslow's hammer states that "If the only tool you have is a hammer, you tend to see every problem as a nail." This is true for most of us when our specialism biases us towards familiar methods and tools. To illustrate this principle, we will look at the problem of gauging the amount of toilet paper that is left on a roll from the outside diameter of that roll and answer the annoying question of why a toilet roll appears to run out so quickly.

We will model the toilet roll with four variables, the external diameter, $d$, of the inner cardboard core, the external diameter, $D$, of the complete roll, the thickness, $t$, of the paper, and the length, $L$, of the remaining paper. If statistical modelling via regression was our only tool, we might start with the equation, $L = f(D, d, t)$, and collect data that would help us to determine the form of the function, $f$, using some form of model fitting. If we use some physics, however, we can simplify the equation into one with fewer variables and this is what the Pi Theorem enables.

The Pi Theorem, which is fundamental to the field of dimensional analysis, is probably unfamiliar to people who work outside of fluid mechanics and engineering, but it has implications for any study that involves physical variables, i.e., variables that are measured in physical (SI) units (The importance of dimensional analysis in statistics has long been understood if not widely applied[34]). In simple terms the Pi Theorem states that a physically meaningful relationship between $n$ physical variables that involves $k$ physical dimensions (i.e. SI base units) can be expressed in terms of $n - k$ dimensionless variables. A full treatment of the Pi Theorem is outside the scope of this article, but it has been covered in many places and is easy to research[35]. For the toilet roll problem the only physical dimension is that of length, so we can reduce the original four variables to three dimensionless quantities. The reduction is not unique but for this problem it is easy to do by inspection. It turns out that a good choice is $\frac{L}{D-d}$, $\frac{D}{t}$, and $\frac{d}{t}$. This gives the new form, $\frac{L}{D-d} = g(\frac{D}{t}, \frac{d}{t})$, for the physical relationship, which reduces the problem from finding a function of three variables to now finding a function in only two variables.

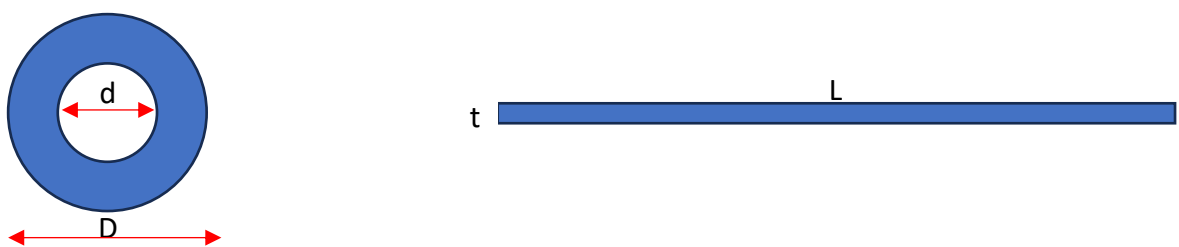

We can deduce the function using simple geometric reasoning. Viewed end-on the area occupied by the toilet paper is $\frac{\pi}{4}(D^2 - d^2)$, when it is rolled up, and, if we visualise it laid flat, the area is $Lt$. Equating the two areas gives $Lt = \frac{\pi}{4}(D^2 - d^2)$, which can be rearranged to give:-

---

[34] For example see Finney (1977) "Dimensions of statistics" *Applied Statistics*, " I am surprised by the lack of attention given to dimensions as a check on the theory and practice of statistics. The basic ideas, readily appreciated, should form part of the stock-in-trade of every statistician."

[35] For example the Wikipedia article is very informative, see  https://en.wikipedia.org/wiki/Buckingham_%CF%80_theorem

Appendix 2 (cont:)

$$\frac{L}{D-d} = \frac{\pi}{4}\left(\frac{D}{t} - \frac{d}{t}\right) = g\left(\frac{D}{t}, \frac{d}{t}\right)$$

which is in the form that we deduced by applying the pi theorem.

Finally, we note that $\frac{dD}{dL} = \frac{2}{\pi}\frac{t}{D}$, which explains why toilet rolls shrink in diameter more quickly as the amount of paper runs out.

# Appendix 3: measuring tree volume.

It is not possible to measure the volume of a tree without felling it so the process begins by felling the tree and removing the crown to leave the stem of the tree.

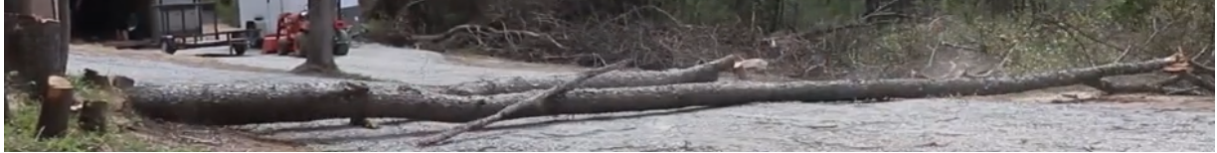

The tree stem is then divided into logs of a convenient length for the mill that will be used to process them into boards.

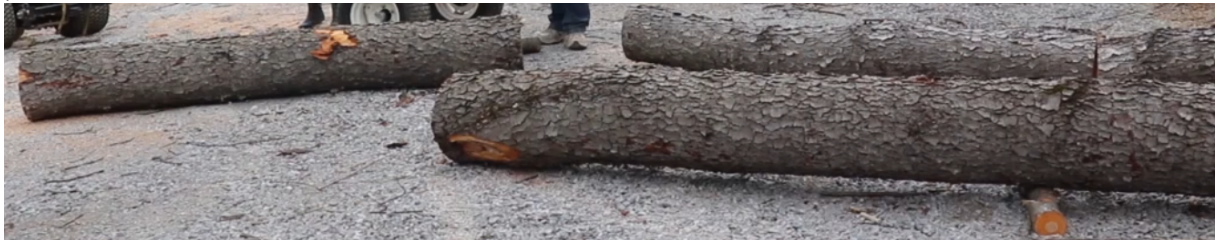

The tree's volume is calculated by measuring the volume of each log and then summing the individual values to get a total value.

The volume of the individual logs is usually calculated using Smalian's formula,
$$V = \frac{\pi}{8}(d_{se}^2 + d_{le}^2)L,$$
where $d_{se}$ is the diameter of the small end, $d_{le}$ is the diameter of the large end, and $L$ is the length of the log.

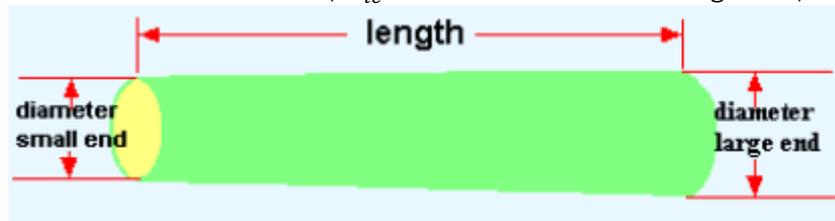